\RequirePackage{amsmath} 
\documentclass[12pt]{iopart}

\usepackage{iopams}  
\usepackage{bbm}
\usepackage{bm}
\usepackage{amsthm}
\usepackage{multirow}
\usepackage{tikz}
\usepackage{pgfplots}
\pgfplotsset{compat=newest}

\newcommand{\sobolev}{H}
\newcommand{\lebesgue}{L}

\renewcommand{\phi}{\varphi}

\newtheorem{thm}{Theorem}

\newtheorem{lem}[thm]{Lemma}
\newtheorem{Cor}[thm]{Corollary}

\begin{document}
\title[Elastic interior transmission eigenvalues]{Elastic interior transmission eigenvalues and their computation via the method of fundamental solutions}

\author{Andreas Kleefeld and
	 Lukas Pieronek\footnote{Author to whom any correspondence should be addressed.}}

\address{Forschungszentrum J\"ulich GmbH, 
	J\"ulich Supercomputing Centre, 52425 J\"ulich, Germany}
\eads{\mailto{a.kleefeld@fz-juelich.de} and \mailto{l.pieronek@fz-juelich.de}}
\vspace{10pt}
\begin{indented}
\item	31 January 2019, revised 4 June 2019
\end{indented}
\begin{abstract}
	A stabilized version of the fundamental solution method to catch ill-conditioning effects is investigated with focus on the computation of complex-valued
	elastic interior transmission eigenvalues in two dimensions for homogeneous and isotropic media without voids. Its algorithm can be implemented very shortly and adopts to many similar PDE-based eigenproblems as long as the underlying fundamental solution function can be easily generated. We develop a corroborative approximation analysis which also implicates new basic results for transmission eigenfunctions and present some numerical examples which together prove successful feasibility of our eigenvalue recovery approach. 
	\\[0.5cm]
\end{abstract}


\section{Introduction}
Interior transmission eigenvalues (ITEs) arise primarily in the study of inverse scattering theory, cf. \cite{bellis2012nature,bellis2010existence,charalambopoulos2006factorization}. Their major role comes along with conventional reconstruction algorithms for which incident wave frequencies of that specific order need to be excluded to fully justify the feasibility of qualitative methods for the target recovery process. While the acoustic and electromagnetic cases have been thoroughly investigated in the last years, see \cite{Colton2012,Sun2016,kleefeld2013numerical,Cakoni2016,kirsch2013inside,kleefeld2018method,kleefeld2018computing,lakshtanov2012bounds}, there are comparably only a few papers covering the numerical study of elastic ITEs so far, e.g. \cite{ji2018computation, peters2017inside,xi2018lowest,xi2018c0ip}, especially when focusing on the overall complex-valued spectrum. It is the purpose of this paper to introduce the approved algorithm from \cite{kleefeld2018computing} and \cite{kleefeld2018method} also for their computation.
As for the latter, we confine to ideally isotropic, homogeneous and planar-symmetric scatterers to resemble the 2D case, but allow for less regular shapes within our theoretical foundation now which are assumed to be known in advance. 

Numerical challenges in computing ITEs origin from the fact that the underlying interior transmission eigenvalue problem (ITP) is non-self-adjoint, 
non-elliptic and non-linear in the eigenvalue parameter. Its straightforward discretization would therefore result in non-Hermitian matrices whose pseudospectra are generally harder to capture, especially when the matrix size becomes large. In this context, the advantage of Ritz-type methods such as the method of fundamental solution (MFS) whose prototype goes back to Kupradze in the 1960th for approximating solutions of well-posed boundary value problems, see \cite{kupradze1965potential}, is that surprisingly accurate results can already be obtained for many regular domains while only a relatively small number of trial functions is used. Due to this observation our MFS-based eigenvalue algorithm similar to \cite{betcke2005reviving} will be very efficient at remarkably low numerical costs and is still easy-to-implement as being mesh- and integration-free unlike the usual competitive methods, cf. \cite{ji2018computation}.

In the course of substantiating our approach mathematically, we prove some novel findings concerning norm relations of ITP eigenfunction pairs which surprisingly differ for purely real and complex-valued ITEs with non-vanishing imaginary part. Since the field of eigenfunction properties is quite undiscovered unlike their eigenvalues themselves in the ITP context, see \cite{blaasten2017vanishing}, we want to point out that most of our results presented here should be adaptable for the acoustic and electromagnetic case, too.

The remainder of this paper is structured as follows: In section 2 we will setup the modeling framework of the ITP. Section 3 will then consider the eigenvalue approximation problem from an abstract and thus more general perspective that also fits other boundary control techniques such as the boundary element method. Definitions and prerequisite results that are specific for our MFS ansatz will be given in section 4 and are followed by numerical examples in section 5. A final conclusion will summarize the merits of our proposed method at the end. 
\section{Setup of the elastic interior transmission eigenvalue problem}
We model time harmonic vibrations of a bounded and elastically penetrable solid $D$ with homogeneous mass density $\rho=\text{const}>0$ which is assumed to differ from its normalized background material through $\rho\ne 1$. Mathematically, we think of $D$ as a domain $D\subset \mathbb{R}^2$ of class $C^{1,1}$ whose regional displacement, when absent or present in comparison with its surrounding, will be modeled by vector fields $u,v: D\to \mathbb{C}^2$, respectively. For both the generalized but isotropic version of Hooke's law shall  relate their gradients in a linear way to the internal stress $\sigma:\mathbb{C}^{2\times2}\to \mathbb{C}^{2\times2}$ via 
\begin{align*}
\sigma(z)
=2\mu\epsilon(z)+\lambda \text{tr}(\epsilon(z))\mathrm{I}\ .
\end{align*}
Here, $\mathrm{I}\in\mathbb{R}^{2\times2}$ is the identity matrix and $\epsilon(z)=(\nabla z+(\nabla z)^{\mathrm{T}})/2$ is the symmetric part of the gradient for any displacement field $z$ measuring its corresponding strain. 
Further, $\mu,\lambda$ are Lam\'{e} parameters which coincide for simplicity for both the scatterer and the background and which are constrained to $\mu>0$ as well as $2\mu+\lambda>0$ to guarantee strong ellipticity of the governing Navier system, see \cite{McLean2000}
\begin{align}\label{navier}
\mathrm{div}(\sigma(z))+\varrho \omega^2 z=\big( \mu \Delta z + (\lambda+\mu)\nabla(\mathrm{div}\,z)\big)+\varrho \omega^2 z=:\Delta^*z+\varrho \omega^2 z=0 \ .
\end{align}
These equations describe from a physical perspective the spatial part of the elastic wave propagation, where $\varrho$ represents any mass density under consideration. 
With this notation at hand we may now formulate the elastic interior transmission problem (ITP):
\begin{align}\label{ITP}
\begin{split}
\Delta^*u+\omega^2 u=0 \quad &\text{in }D\ , \\
\Delta^*v+\rho\omega^2 v=0 \quad &\text{in }D\ , \\
u=v \quad &\text{on }\partial D\ , \\
\sigma(u)\nu=\sigma(v)\nu  \quad &\text{on }\partial D\ ,
\end{split}
\end{align}
where $u,v\in \lebesgue^2(D,\mathbb{C}^2)$ and the scattered part $(u-v)$ is supposed to be an element of the more regular Sobolev space
\begin{align*}
H^2_0(D)=\{\phi\in \sobolev^2(D,\mathbb{C}^2): \phi=0 \text{ and } \sigma(\phi)\nu=0 \text{ on }\partial D \}
\end{align*}
with outer normal $\nu$ along the boundary $\partial D$. Frequencies $\omega\in \mathbb{C}\backslash \{0\}$ which admit non-trivial solutions $(u,v)$ to the above ITP will be called interior transmission eigenvalues (ITEs). Apparently,
they refer to those pairs of harmonic waves in $D$ whose behavior along the scattering boundary coincides and which thus lock the possibility of detecting the scatterer on the basis of close-by data.

Note that the co-normal derivative appearing in the Neumann boundary condition of \eqref{ITP} is connected naturally to the highest order term of their PDEs via integration by parts, also known as Betti's first formula
\begin{align}
\begin{split}\label{betti1}
\int_D\Delta^*\phi\cdot\psi\ \mathrm{d}x&=-\int_D\sigma(\phi) : \nabla \psi\ \mathrm{d}x+\int_{\partial D}(\sigma(\phi)\nu)\cdot\psi\ \mathrm{d}s\\
&= -\int_D 2\mu \epsilon(\phi):\epsilon(\psi)+\lambda\mathrm{div}\,\phi\,\mathrm{div}\,\psi\ \mathrm{d}x+\int_{\partial D}(\sigma(\phi)\nu)\cdot\psi\ \mathrm{d}s
\end{split}
\end{align}
which holds by duality for any $\psi\in L^2(D,\mathbb{C}^2)$ and $\phi\in \sobolev^2(D,\mathbb{C}^2)$. Here, the colon symbol denotes the Frobenius inner product given by $A:B=\mathrm{tr}(AB^{\mathrm{H}})$, whereas the single dot refers to the scalar-product-like operation $a\cdot b=a^{\mathrm{H}}b$ for $a,b\in \mathbb{C}^2$. Having thus set the mathematical framework of our eigenproblem, in the sequel we try to recover ITEs $\omega$ by solving \eqref{ITP} approximately in the sense that we allow for small deviations within the boundary data that are assumed to vanish in some limiting procedure.

\section{Approximation analysis of elastic interior transmission eigenvalues via boundary control}
We define our relaxed space of trial functions for approximating solutions of \eqref{ITP} subject to boundary optimization by
\begin{align}\label{ansatz}
\mathcal{H}:=\bigcup_{0\leq\arg(\omega)< \frac{\pi}{4}}\mathcal{H}(\omega)\ ,
\end{align}
where for any $\omega\in \mathbb{C}\backslash \{0\}$
\begin{align*}
\mathcal{H}(\omega):= \big\{(u,v)\in C^{\infty}(\overline{D})\times C^{\infty}(\overline{D}):\ \Delta^*u+\omega^2 u=0\ ,\  \Delta^*v+\omega^2 \rho v=0\big\}\ .
\end{align*}
By definition, any pair in $\mathcal{H}(\omega)$ fulfills the required PDE conditions from \eqref{ITP} automatically and it is the choice of $\omega$ that determines in how far their boundary data are compatible in the sense of the ITP.
Note that the above restrictions on $\omega$ in $\mathcal{H}$ are due to the following theorem on the overall locations of complex ITEs and the fact that both $\overline{\omega}$, $-\overline{\omega}$ and $-\omega$ are each ITEs if and only if some $\omega$ from the first quadrant in the complex plane is. We drop the proof of the former since it would exactly follow the lines from its acoustic analogon in \cite{cakoni2010interior} with the obvious operator adaptions.  
\begin{thm}
	Let $\omega=\omega_1+\mathrm{i}\omega_2\in \mathbb{C}\backslash \{0\}$	be an ITE with $\omega_1,\omega_2\in \mathbb{R}$ for the scatterer $D$. Then it holds that 
	\begin{align*}
	\omega_1^2>\omega_2^2 \qquad \text{and} \qquad \omega_1^4+\omega_2^4+\frac{2\rho+6}{\delta-1}\omega_1\omega_2-\frac{\lambda(D)}{\rho}(\omega_1^2-\omega_2^2)>0 \ ,
	\end{align*}
	where $\lambda(D)$ is the smallest interior Dirichlet eigenvalue for the Navier problem \eqref{navier} with $\varrho=\rho$. In particular, if $\omega$ lies in the first quadrant of the complex plane, then $0\leq\arg(\omega)< \pi/4$.
\end{thm}

We aim to extract those $\omega$ which allow for approximate ITP eigenfunctions in $\mathcal{H}$ with a relatively small ratio of boundary misfit to interior norm. The next theorem states that if these residual quotients can be made arbitrarily small while the corresponding $\omega$ accumulate, their limit is indeed an ITE. Its proof slightly refines the technique from \cite[Theorem 2]{kleefeld2018method} to encompass also complex-valued eigenvalues now.
\begin{thm}\label{theorem1}
	Assume that $\{(u_m,v_m,\omega_m)\}_{m \in \mathbb{N}}\subset \mathcal{H}\times \mathbb{C}$ fulfill for some $1\leq C<\infty$ the following conditions:
	\begin{enumerate}	
		\item 	eigenvalue convergence: $\omega_m \to \omega$ such that $\mathrm{arg}(\omega)<1/4$,
		
		\item	uniform interior bound: $1/C\le\left(\|u_m\|^2_{\lebesgue^2(D,\mathbb{C}^2)}+\|v_m\|^2_{\lebesgue^2(D,\mathbb{C}^2)}\right)\le C$ for all $m$ large enough,
		
		\item	vanishing boundary misfit: 
		$\Big(\|u_m-v_m\|_{\sobolev^{\frac{3}{2}}(\partial D,\mathbb{C}^2)}+\|\sigma(u_m-v_m)\nu\|_{\sobolev^{\frac{1}{2}}(\partial D,\mathbb{C}^2)}\Big)\to 0$ when $m\to \infty$.
	\end{enumerate}
	Then, $\omega$ is an ITE and a subsequence of $(u_m,v_m)$ converges weakly to some eigenfunction pair.
	\begin{proof}
		By rescaling and redefining $(u_m,v_m)$, if necessary, we can assume without loss of generality that $C=1$ and aim to apply weak compactness in order to construct a solution candidate which will indeed meet all the required ITE criteria. 
		By assumption $(ii)$ we know (modulo the extraction of subsequences which we will not relabel in $m$) that $u_m\rightharpoonup u$ and $v_m\rightharpoonup v$ in $L^2(D,\mathbb{C}^2)$ which implies that $(u,v)\in L^2(D,\mathbb{C}^2)\times L^2(D,\mathbb{C}^2)$ fulfills the interior conditions of the ITP, or equivalently the PDE system \eqref{navier}, according to
		\begin{align*}
		\int_D u\cdot(\Delta^*\phi + \omega^2 \phi)\, \mathrm{d}x&=\lim_{m\to \infty}\int_D u_m\cdot(\Delta^*\phi + \omega_m^2 \phi)\, \mathrm{d}x \\
		&=\lim_{m\to \infty}\int_D (\Delta^*u_m + \omega^2 u_m)\cdot\phi\,\mathrm{d}x=0
		\end{align*}
		for any bump function $\phi\in C_c^{\infty}(D,\mathbb{C}^2)$ and with a similar calculus for $v$. In order to prove that $(u-v)$ also has the correct ITP boundary data, it suffices to prove that the differences $(u_m-v_m)$ are bounded in $\sobolev^2(D,\mathbb{C}^2)$ in combination with assumption $(iii)$ because $w\mapsto w_{|\partial D}$ and $w\mapsto \sigma(w)_{|\partial D}\nu$ are continuous as operators $\sobolev^2(D,\mathbb{C}^2)\to \sobolev^{\frac{3}{2}}(\partial D,\mathbb{C}^2)$ and $\sobolev^2(D,\mathbb{C}^2)\to \sobolev^{\frac{1}{2}}(\partial D,\mathbb{C}^2)$, respectively, due to our regularity assumptions on $\partial D$. 
		Noting that $(u_m-v_m)$ are solutions to the inhomogeneous Navier system
		\begin{align*}
		\Delta^*(u_m-v_m) + \omega_m^2 (u_m-v_m)=(1-\rho)v_m \quad \text{in }D\ ,
		\end{align*}
		elliptic estimates like in \cite{McLean2000} tell us that 
		\begin{align}\label{homogeneous}
		\|u_m-v_m\|_{\sobolev^2(D,\mathbb{C}^2)}\leq C \left(\|u_m-v_m\|_{\sobolev^{\frac{3}{2}}(\partial D,\mathbb{C}^2)}+\|u_m\|_{\lebesgue^2(D,\mathbb{C}^2)}+\|v_m\|_{\lebesgue^2(D,\mathbb{C}^2)}\right)
		\end{align}
		which therefore gives the desired uniform bound with respect to $m$.
		
		It remains to show that $(u,v)$ is non-trivial. For this we recall that the embedding $\sobolev^2(D,\mathbb{C}^2)\hookrightarrow \lebesgue^2(D,\mathbb{C}^2)$ is compact which implies $(u_m-v_m)\to (u-v)$ strongly in $\lebesgue^2(D,\mathbb{C}^2)$. Apparently, $(u,v)\ne 0$ if $\|u-v\|_{\lebesgue^2(D,\mathbb{C}^2)}>0$ so we will assume contrarily that $(u_m-v_m)\to 0$ in $\lebesgue^2(D,\mathbb{C}^2)$. 
		Then, on the one hand, the bounded sequence
		\begin{align*}
		a_m:=\int_Du_m\cdot\overline{v}_m\, \mathrm{d}x
		\end{align*}
		may be singled out to converge to some $a\in \mathbb{C}$ from which we know by $(ii)$ with $C=1$ that
		\begin{align}\label{plot}
		\frac{1}{2}\geq|a|\geq\mathrm{Re}\,a=\lim_{m\to \infty}\frac{\|u_m\|^2_{\lebesgue^2(D,\mathbb{C}^2)}+\|v_m\|^2_{\lebesgue^2(D,\mathbb{C}^2)}-\|u_m-v_m\|^2_{\lebesgue^2(D,\mathbb{C}^2)}}{2}\geq \frac{1}{2}\ ,
		\end{align}
		i.e. $a=1/2$. Since $\mathrm{arg}(\omega)<1/4$ and $\rho\ne 1$, we can even conclude $|\mathrm{Re}\,a(\omega^2-\rho\overline{\omega}^2)|>0$.

		On the other hand, \eqref{betti1} with exchanged roles of its test functions $\phi,\psi \in C^{\infty}(\overline{D})$ yields the analogon of Green's second identity
		\begin{align*}
		\int_D\psi\cdot\Delta^*\phi-\phi\cdot\Delta^*\psi\,\mathrm{d}x=\int_{\partial D}\psi\cdot(\sigma(\phi)\nu)-\phi\cdot(\sigma(\psi)\nu)\,\mathrm{d}s\ .
		\end{align*}
		Duality shows that the separated boundary contributions of the approximate eigenfunctions $\|u_m\|_{\sobolev^{-\frac{1}{2}}(\partial D,\mathbb{C}^2)}$,$\|v_m\|_{\sobolev^{-\frac{1}{2}}(\partial D,\mathbb{C}^2)}$,$\|\sigma(u_m)\nu\|_{\sobolev^{-\frac{3}{2}}(\partial D,\mathbb{C}^2)}$ and $\|\sigma(v_m)\nu\|_{\sobolev^{-\frac{3}{2}}(\partial D,\mathbb{C}^2)}$ are even uniformly bounded in $m$ thanks to $(i)$--$(iii)$. Therefore we may further compute
		\begin{align*}
		0&<|\mathrm{Re}\,a(\omega^2-\rho\overline{\omega}^2)| \\
		&=\lim_{m\to \infty}\left|\mathrm{Re}\int_D(\omega^2-\rho\overline{\omega}^2)u_m\cdot\overline{v}_m\,\mathrm{d}x\right|\\
		&=\lim_{m\to \infty}\left|\mathrm{Re}\int_D(\omega_m^2u_m)\cdot\overline{v}_m-u_m\cdot(\rho\overline{\omega}_m^2\overline{v}_m)
		\,\mathrm{d}x\right|\\
		&=\lim_{m\to \infty}\left|\mathrm{Re}\int_D\overline{v}_m\cdot\Delta^* u_m-u_m\cdot\Delta^*\overline{v}_m
		\,\mathrm{d}x\right| \\
		&=\lim_{m\to \infty}\left|\mathrm{Re}\int_{\partial D}\overline{v}_m\cdot(\sigma(u_m)\nu)-u_m\cdot(\sigma(\overline{v}_m)\nu)
		\,\mathrm{d}s\right| \\
		&=\lim_{m\to \infty}\left|\mathrm{Re}\int_{\partial D}\overline{v}_m\cdot(\sigma(u_m)\nu)-\overline{u}_m\cdot(\sigma(v_m)\nu)
		\,\mathrm{d}s\right| \\
		&=\lim_{m\to \infty}\left|\mathrm{Re}\int_{\partial D}(\overline{v}_m-\overline{u}_m)\cdot(\sigma(u_m)\nu)   
		+\overline{u}_m\cdot(\sigma(u_m-v_m)\nu)\,\mathrm{d}s\right| \\
		&\leq\lim_{m\to \infty}
		\|v_m-u_m\|_{\sobolev^{\frac{3}{2}}(\partial D,\mathbb{C}^2)}
		\|\sigma(u_m)\nu\|_{\sobolev^{-\frac{3}{2}}(\partial D,\mathbb{C}^2)}\\
		& \ \ \ +\lim_{m\to \infty}
		\|\sigma(v_m-u_m)\nu\|_{\sobolev^{\frac{1}{2}}(\partial D,\mathbb{C}^2)}
		\|u_m\|_{\sobolev^{-\frac{1}{2}}(\partial D,\mathbb{C}^2)}\\
		&=0\ ,
		\end{align*}
		which is a contradiction.
	\end{proof}
\end{thm}
In order to establish for fixed $m$ some more qualitative relation between the boundary misfit of some sufficiently approximate eigenfunction pair and its absolute eigenvalue deviation from the actual one, we derive an interconnecting estimate in the following which can also be seen as practical guiding principle for numerical calculations. 
Similar to the acoustic case, see \cite[Lemma 5]{kleefeld2018method}, but now improved for even less regular shapes, its rigorous validity relies on a critical integral expression which must not vanish. 
\begin{lem}\label{lemma}
	Let $(u,v)$ be an ITP eigenfunction pair with ITE $\omega$ and assume that $(\widetilde{u},\widetilde{v})\in \mathcal{H}(\widetilde{\omega})$ for some aritrary frequency $w$. 
	If the integral constraint 
	\begin{align}\label{constraint}
	\Bigg|\int_D u\cdot\widetilde{u}-\rho v\cdot\widetilde{v}\, \mathrm{d}x\Bigg|\geq \widetilde{\varepsilon}> 0
	\end{align}
	is fulfilled, then there exists a constant $C>0$ determined only by the boundary data of $u$ (or equivalently $v$) such that
	\begin{align}\label{error}
	|\omega^2-\widetilde{\omega}^2|\leq \frac{C}{\widetilde{\varepsilon}} \sqrt{\|\widetilde{u}-\widetilde{v}\|^2_{\sobolev^{\frac{3}{2}}(\partial D,\mathbb{C}^2)}+\|\sigma(\widetilde{u}-\widetilde{v})\nu\|^2_{\sobolev^{\frac{1}{2}}(\partial D,\mathbb{C}^2)}}\ . 
	\end{align}
	\begin{proof}
		Applying Betti's formula twice, using the identical ITP boundary data for $u$ and $v$, yields
		\begin{align*}
		&\omega^2\int_{D}\widetilde{u}\cdot u-\rho\widetilde{v}\cdot v\, \mathrm{d}x\\
		=& \int_{D}-\widetilde{u}\cdot\Delta^* u+\widetilde{v}\cdot\Delta^* v\, \mathrm{d}x\\
		=& \int_D 2\mu \epsilon(\widetilde{u}):\epsilon(u)+\lambda\mathrm{div}\,\widetilde{u}\,\mathrm{div}\,u\ \mathrm{d}x
		-\int_D 2\mu \epsilon(\widetilde{v}):\epsilon(v)+\lambda\mathrm{div}\,\widetilde{v}\,\mathrm{div}\,v\ \mathrm{d}x \\ &-\int_{\partial D}(\widetilde{u}-\widetilde{v})\cdot(\sigma (u)\nu)\, \mathrm{d}s\\
		=& \int_{D}(\Delta^* \widetilde{u}\cdot u+\Delta^* \widetilde{v}\cdot v)\, \mathrm{d}x
		+\int_{\partial D}(\sigma(\widetilde{u}-\widetilde{v})\nu)\cdot u\, \mathrm{d}s
		-\int_{\partial D}(\widetilde{u}-\widetilde{v})\cdot(\sigma(u)\nu)\, \mathrm{d}s\\
		=& \ \widetilde{\omega}^2\int_{D}\widetilde{u}\cdot u-\rho\widetilde{v}\cdot v\, \mathrm{d}x
		+\int_{\partial D}(\sigma(\widetilde{u}-\widetilde{v})\nu)\cdot u\, \mathrm{d}s
		-\int_{\partial D}(\widetilde{u}-\widetilde{v})\cdot(\sigma(u)\nu)\, \mathrm{d}s
		\end{align*}
		and after rearranging we obtain
		\begin{align*}
		(\omega^2-\widetilde{\omega}^2)\int_{D}\widetilde{u}\cdot u-\rho\widetilde{v}\cdot v\, \mathrm{d}x=\int_{\partial D}(\sigma(\widetilde{u}-\widetilde{v})\nu)\cdot u\, \mathrm{d}s
		-\int_{\partial D}(\widetilde{u}-\widetilde{v})\cdot(\sigma (u)\nu)\, \mathrm{d}s\ .
		\end{align*}
		Taking absolute values gives 
		\begin{align*}
		&\big|\omega^2-\widetilde{\omega}^2\big|\\
		\leq&\  \frac{1}{\widetilde{\varepsilon}}\left(\int_{\partial D}|(\sigma(\widetilde{u}-\widetilde{v})\nu)\cdot u|\, \mathrm{d}s+
		\int_{\partial D}|(\widetilde{u}-\widetilde{v})\cdot(\sigma (u)\nu)|\, \mathrm{d}s\right)\\
		\leq&\  \frac{1}{\widetilde{\varepsilon}}\left(\|\sigma(\widetilde{u}-\widetilde{v})\nu\|_{\sobolev^{\frac{1}{2}}(\partial D,\mathbb{C}^2)}\|u\|_{\sobolev^{-\frac{1}{2}}(\partial D,\mathbb{C}^2)}+\|\widetilde{u}-\widetilde{v}\|_{\sobolev^{\frac{3}{2}}(\partial D,\mathbb{C}^2)}\|\sigma (u)\nu\|_{\sobolev^{-\frac{3}{2}}(\partial D,\mathbb{C}^2)}\right)\\
		\leq&\  \frac{C}{\widetilde{\varepsilon}} \sqrt{\|\widetilde{u}-\widetilde{v}\|_{\sobolev^{\frac{3}{2}}(\partial D,\mathbb{C}^2)}^2
			+\|\sigma(\widetilde{u}-\widetilde{v})\nu\|_{\sobolev^{\frac{1}{2}}(\partial D,\mathbb{C}^2)}\|^2}\ ,
		\end{align*}
		where
		\begin{align*}
		C:=\sqrt{\|u\|_{\sobolev^{-\frac{1}{2}}(\partial D,\mathbb{C}^2)}^2+\|\sigma (u)\nu\|_{\sobolev^{-\frac{3}{2}}(\partial D,\mathbb{C}^2)}^2}\ .
		\end{align*}
		Duality shows that $C<\infty$ since $u$ solves the Navier equation.
	\end{proof}
\end{lem}

We state a direct consequence with respect to the $L^2(D,\mathbb{C}^2)$-norm of eigenfunctions for frequencies $\omega\in \mathbb{C}\backslash \mathbb{R}$.

\begin{Cor}\label{Coro}
	If $(u,v)$ is an ITP eigenfunction pair with $\omega\in \mathbb{C}\backslash \mathbb{R}$, then we have that $\|u\|^2_{L^2(D,\mathbb{C}^2)}=\rho\|v\|^2_{L^2(D,\mathbb{C}^2)}$.
	\begin{proof}
		Let $\{(u_m,v_m)\}_{m\in \mathbb{N}}$ be a sequence in $\mathcal{H}(\omega)$ such that $w_m:=(u_m-v_m)\to(u-v)$ in $H^2(D,\mathbb{C}^2)$, see for example Theorem \ref{density}\ later. In particular, $(\overline{u}_m,\overline{v}_m)\in\mathcal{H}(\overline{\omega})$, $\|\overline{w}_m\|^2_{\sobolev^{\frac{3}{2}}(\partial D,\mathbb{C}^2)}+\|\sigma(\overline{w}_m)\|^2_{\sobolev^{\frac{1}{2}}(\partial D,\mathbb{C}^2)}\to 0$ and the left hand side of \eqref{constraint} converges because of $\overline{u}_m=(\Delta^*\overline{w}_m+\rho\overline{\omega}^2\overline{w}_m)/(\overline{\omega}(1-\rho))$ and $\overline{v}_m=(\Delta^*\overline{w}_m+\overline{\omega}^2\overline{w}_m)/(\overline{\omega}(1-\rho))$ to
		\begin{align*}
		\|u\|^2_{L^2(D,\mathbb{C}^2)}-\rho\|v\|^2_{L^2(D,\mathbb{C}^2)}
		\end{align*}
		when setting $\widetilde{u}:=\overline{u}_m$ and $\widetilde{v}:=\overline{v}_m$. Since $|\omega^2-\overline{\omega}^2|>0$ by assumption, the right hand side of \eqref{error} forces any uniform bound on $\widetilde{\varepsilon}$ to vanish in the limit $m\to \infty$.
	\end{proof}
\end{Cor}
The situation is different for $\omega\in \mathbb{R}$ and  especially when it is the ITE of smallest magnitude. In this case, with the additional assumption that $\rho>1$ or $0<\rho<1$ is sufficiently large or small to be made precise next, respectively, we can show that $\|u\|^2_{L^2(D,\mathbb{C}^2)}-\rho\|v\|^2_{L^2(D,\mathbb{C}^2)}\ne 0$ which also guarantees a uniform bound $\widetilde{\varepsilon}>0$ for approximate eigenfunction pairs $(\widetilde{u},\widetilde{v})$ in the vicinity of $(u,v)$ (cf. Corollary \ref{Coro2} later). 
To avoid misleading confusion at that point, we want to emphasize that Corollary \ref{Coro} does not imply the non-existence of some positive constant $\widetilde{\varepsilon}$ when dealing with complex-valued ITEs $\omega$ admitting a non-vanishing imaginary part.
\begin{thm}
	Let $\omega$ be the smallest real-valued ITE for the scatterer $D$ with constant density $\rho\ne 1$. If $\rho>1$ is large or $0<\rho<1$ is small enough, where the corresponding thresholds depend only on $D$ and the Lam\'{e} parameters 
	$\mu, \lambda$, we have the relations
	\begin{align*}
	\|u\|^2_{L^2(D,\mathbb{C}^2)}-\rho\|v\|^2_{L^2(D,\mathbb{C}^2)}<0 \qquad \text{or} \qquad \|u\|^2_{L^2(D,\mathbb{C}^2)}-\rho\|v\|^2_{L^2(D,\mathbb{C}^2)}>0\ ,
	\end{align*}
	respectively.
	\begin{proof}
		For the sake of presentation we will assume that $\rho>1$ since the case $0<\rho<1$ works structurally similar. Because $u$ and $v$ can be expressed each in terms of their difference $w:=u-v$ by $u=(\Delta^*w+\rho\omega^2w)/(\omega^2(\rho-1))$ and $v=(\Delta^*w+\omega^2w)/(\omega^2(\rho-1))$, the basic idea of our proof will be to exploit isometry of the Fourier transform with respect to the single field $w$ to obtain a well-behaved algebraic integrand in terms of $\rho$.
		The fact $w\in H^2_0(D)$ then shows that $u,v,w$ extend naturally by zero outside of $D$ so Plancherel's identity gives
		\begin{align*}
		\int_D |u|^2-\rho|v|^2\,\mathrm{d}x
		&=\frac{1}{(2\pi)^2}\int_{\mathbb{R}^2} |\mathcal{F}u|^2-\rho|\mathcal{F}v|^2\,\mathrm{d}\xi \\
		&=\frac{1}{(2\pi)^2}\int_{\mathbb{R}^2}\frac{|\mathcal{F}(\Delta^*w+\rho\omega^2w)|^2-\rho |\mathcal{F}(\Delta^*w+\omega^2w)|^2}{(\rho-1)^2\omega^4}\,\mathrm{d}\xi\\
		&=\frac{1}{(2\pi)^2}\int_{\mathbb{R}^2}\frac{-|\mathcal{F}\Delta^*w|^2+\rho\omega^4|\mathcal{F}w|^2}{(\rho-1)\omega^4}\,\mathrm{d}\xi\\
		&\leq \frac{1}{(2\pi)^2}\int_{\mathbb{R}^2} \frac{-c|\xi|^4+\rho\omega^4}{(\rho-1)\omega^4}|\mathcal{F}w|^2\,\mathrm{d}\xi\ .
		\end{align*}
		In the last step we employed the pointwise estimate
		\begin{align*}
		|\mathcal{F}\Delta^*w(\xi)|^2\geq c |\xi|^4|\mathcal{F}w(\xi)|^2
		\end{align*}
		with $c:=\min\{\mu,\lambda+2\mu\}>0$ inherited from strong ellipticity. As will be shown later, we have for the smallest ITE $\omega=\omega({\rho})$ with density $\rho$ that $\rho\omega^4\to 0$ and likewise $R(\rho)\to 0$ for $\rho\to \infty$, where $R(\rho):=\omega\sqrt[4]{(2\rho-1)/c}$ is the solution of $p(t):=(-ct^4+\rho\omega^4)/(\rho-1)\omega^4=-1$ in $t$. With $B_{R(\rho)}$ being the disc centered at the origin of radius $R(\rho)$, we split the integral above and exploit monotonicity of $p$ to deduce
		\begin{align*}
		&\int_{B_{R(\rho)}} \frac{-c|\xi|^4+\rho\omega^4}{(\rho-1)\omega^4}|\mathcal{F}w|^2\,\mathrm{d}\xi+\int_{\mathbb{R}^2\backslash B_{R(\rho)}} \frac{-c|\xi|^4+\rho\omega^4}{(\rho-1)\omega^4}|\mathcal{F}w|^2\,\mathrm{d}\xi \\
		\leq&\ p(0)\int_{B_{R(\rho)}} |\mathcal{F}w|^2\,\mathrm{d}\xi+p(R(\rho))\int_{\mathbb{R}^2\backslash B_{R(\rho)}} |\mathcal{F}w|^2\,\mathrm{d}\xi \\
		=& \frac{\rho}{\rho-1}\|\mathcal{F}w\|^2_{L^2(B_{R(\rho)},\mathbb{C}^2)}-\left(\|\mathcal{F}w\|^2_{L^2(\mathbb{R}^2,\mathbb{C}^2)}-\|\mathcal{F}w\|^2_{L^2(B_{R(\rho)},\mathbb{C}^2)}\right) \\
		\leq& \left(\frac{\rho}{\rho-1}+1\right)\|\mathcal{F}w\|^2_{L^2(B_{R(\rho)},\mathbb{C}^2)}-(2\pi)^2\|w\|^2_{L^2(D,\mathbb{C}^2)}\ .
		\end{align*}
		The first summand can be made arbitrarily small in terms of $\rho$ because of
		\begin{align*}
		\|\mathcal{F}w\|^2_{L^2(B_{R(\rho)},\mathbb{C}^2)}
		&\leq \big(\max |\mathcal{F}w|\big)^2\pi R(\rho)^2\\
		&\leq \|w\|^2_{L^1(D,\mathbb{C}^2)}\pi R(\rho)^2\\
		&\leq  \|w\|^2_{L^2(D,\mathbb{C}^2)}R(\rho)^2\mathcal{L}^2(D)\ ,
		\end{align*}
		where $\mathcal{L}^2(D)$ denotes the two-dimensional Lebesgue measure of $D$. Putting everything together, we finally obtain
		\begin{align*}
		\int_D |u|^2-\rho|v|^2\,\mathrm{d}x
		=\|w\|^2_{L^2(D,\mathbb{C}^2)}\left(\left(\frac{\rho}{\rho-1}+1\right)\mathcal{L}^2(D)R(\rho)^2-(2\pi)^2\right)<0
		\end{align*}
		for $\rho$ large enough due to the decay of $R(\rho)$.
		
		As announced for the latter, it remains to show that $\rho\omega^4\to 0$ as $\rho\to \infty$, where $\omega$ is the smallest real-valued ITE of $D$ with density $\rho$. According to Corollary 1 from \cite{bellis2012nature}, the magnitude of $\omega$ can be bounded from above by the smallest ITE from any included disc $B_r\subset D$ which thus amounts to show the asymptotics assertion just for the unit disc as scatterer. In this case, as was derived in the appendix of \cite{ji2018computation} with an ansatz for a purely compressional wave instead of its orthogonal shear part that we will consider now for completion, $\omega$ is some ITE corresponding to a radial symmetric eigenfunction if and only if 
		\begin{align}\label{det}
		\mathrm{det}
		\begin{pmatrix}
		J_1\left(\omega \sqrt{\frac{1}{\mu}}\right) & J_1\left(\omega \sqrt{\frac{\rho}{\mu}}\right) \\
		\omega \sqrt{\frac{1}{\mu}}J'_1\left(\omega \sqrt{\frac{1}{\mu}}\right) 
		& \omega \sqrt{\frac{\rho}{\mu}}J'_1\left(\omega \sqrt{\frac{\rho}{\mu}}\right)
		\end{pmatrix}
		=0\ ,
		\end{align} 
		where $J_1$ is the first Bessel function of order one. This condition can be restated as finding roots $\omega$ of the piecewise continuous function
		\begin{align*}
		g(\omega):=f(\omega)-f(\sqrt{\rho}\omega)
		\end{align*}
		which can be recursively decomposed into 
		\begin{align*}
		f(x)=h\left(\frac{x}{\sqrt{\mu}}\right) \qquad \text{and} \qquad 
		h(y):=\frac{yJ_1'(y)}{J_1(y)}\ .
		\end{align*}
		Let $j_1<j_2$ be the two smallest positive roots of $J_1$ and choose $\rho>1$ large enough to have $j_2<j_1\sqrt{\rho}$. Then set $\omega_1:=j_1\sqrt{\mu/\rho}$ as well as $\omega_2:=j_2\sqrt{\mu/\rho}$ and observe that $g$ is singular at those points, but continuous in between. Also, those poles have different signs according to
		\begin{align*}
		\lim_{\omega\searrow \omega_1}g(\omega)
		=-\lim_{\omega\searrow \omega_1}f(\sqrt{\rho}\omega)=-\infty \qquad \text{and} \qquad \lim_{\omega\nearrow \omega_2}g(\omega)
		=-\lim_{\omega\nearrow \omega_2}f(\sqrt{\rho}\omega)=\infty\ ,
		\end{align*} 
		which follows from the basic facts that $J_1<0$ in $(j_1,j_2)$, $J_1'(j_1)<0$ but $J_1'(j_2)>0$ and that both $f(\omega_1)$, $f(\omega_2)$ are finite. Therefore we can make use of the intermediate value theorem which guarantees for any large $\rho$ a root $\omega$ of $g$ 
		fulfilling $\omega_1\leq \omega \leq \omega_2$ or equivalently the uniform bound $j_1^2\mu\leq \rho\omega^2 \leq j_1^2\mu$ as $\rho\to \infty$. In particular, $\rho\omega^4\to 0$ for the same limit procedure which finally proves our lemma.
	\end{proof}

\end{thm}
Finally, we cite a result from \cite[Corollary 6]{kleefeld2018method} that summarizes the previous findings for our conceptual ITE recovery approach. 

\begin{Cor}\label{Coro2}
	Let the conditions of Theorem \ref{theorem1} hold for $\{(u_m,v_m,\omega_m)\}_{m\in\mathbb{N}}\subset \mathcal{H}\times \mathbb{C}$ with ITE $\omega$ and eigenfunction pair $(u,v)\in L^2(D,\mathbb{C}^2)\times L^2(D,\mathbb{C}^2)$. Assume additionally that 
	\begin{align}\label{magnitude}
	\int_D (u^2-\rho v^2)\, \mathrm{d}x\ne 0\ .
	\end{align}
	Then, for sufficiently large $m\in \mathbb{N}$, we have (modulo the relabeling of the weakly convergent subsequence)
	\begin{align*}
	|\omega^2-\omega_m^2|\leq C \sqrt{\|u_m-v_m\|^2_{\sobolev^{\frac{3}{2}}(\partial D)}+\|\sigma(u_m-v_m)\nu\|^2_{\sobolev^{\frac{1}{2}}(\partial D)}}\ ,
	\end{align*}
	where $C>0$ depends on the boundary data of $u$ (or equivalently $v$) and the magnitude of \eqref{magnitude}.
	\begin{proof}
		$(u_m,v_m)\rightharpoonup (u,v)$ in $L^2(D,\mathbb{C}^2)\times L^2(D,\mathbb{C}^2)$ implies with \eqref{magnitude} that $\widetilde{\varepsilon}>0$ in Lemma \ref{lemma} uniformly for $m$ large enough and thus the existence of $C>0$.
	\end{proof} 
\end{Cor}

Altogether, we have proved that the conditioned process from Theorem \ref{theorem1} for detecting ITEs is spurious free as limiting procedure and provided a constrained a posteriori estimate for the eigenvalue approximation accuracy at each step $m$. To meet the conditions therein, we will construct ``simple'' trial functions in the next section that are unspecified in \eqref{ansatz} so far.
\section{The method of fundamental solutions}
\subsection{The abstract setting}
In this section we want to focus on how to generate interior solutions to a given PDE on the basis of its fundamental solution. A fundamental solution
for the free space Navier system \eqref{navier} with arguments $x\ne y \in \mathbb{R}^2$, constant coefficients $\mu,\lambda, \varrho$ and unknown eigenfrequency parameter $\omega$ is given by
\begin{align}\label{fundamental}
&\Phi_{\varrho\omega^2}(x,y)\nonumber \\
:=\, &\Phi_{\varrho\omega^2}^{\mu,\lambda}(|x-y|)\nonumber \\
:=\, &\frac{\mathrm{i}}{4\mu}H^{(1)}_0(k_s|x-y|)\cdot\mathrm{I}+\frac{\mathrm{i}}{4\omega^2}\nabla^{\top} \nabla\left(H^{(1)}_0(k_p|x-y|)-H^{(1)}_0(k_s|x-y|)\right)\ ,
\end{align}
with wave numbers
\begin{align*}
k_s^2:=\frac{\varrho\omega^2}{\mu}\ ,\qquad k_p^2:=\frac{\varrho\omega^2}{\lambda+2\mu}
\end{align*} 
and $H^{(1)}_0$ being the first Hankel function of order zero. 
The indices $s$ and $p$ originate from the well-known Helmholtz decomposition that divides any properly decaying solution $e=e_p+e_s$ of the exterior Navier problem into its compressional and shear wave field. More precisely, $e_p=(-1/k_p^2)\nabla(\mathrm{div}\;e)$ and $e_s=(-1/k_s^2)\big(-\partial_2(\partial_1 e_2-\partial_2 e_1),\partial_1(\partial_1 e_2-\partial_2 e_1)\big)^{\top}$ solve 
\begin{align*}
\Delta e_s +k_s^2e_s=0 \qquad \text{and} \qquad
\Delta e_p +k_p^2e_p=0 \qquad \text{in }\mathbb{R}^2\backslash D
\end{align*}
and are constrained to fulfill Kupradze's (outgoing) radiation condition in two dimensions
\begin{align*}
\lim_{r\to\infty}\sqrt{r}(\partial_re_s-\mathrm{i}k_se_s)=0 \ ,\qquad 
\lim_{r\to\infty}\sqrt{r}(\partial_re_p-\mathrm{i}k_pe_p)=0
\end{align*}
uniformly in all directions of the radial distance $r$ from the origin. In particular, freezing one of the arguments in \eqref{fundamental} such as $y$ without loss of generality, these properties apply column-wise (and by symmetry of the fundamental matrix also row-wise)
to $x\mapsto \Phi_{\varrho\omega^2}(x,y)$ in $\mathbb{R}^2\backslash \{y\}$ for any source point $y\in \mathbb{R}^2$ and density $\varrho=\mathrm{const}$. Choosing some simply closed and sufficiently smooth contour $\Gamma\in \mathbb{R}^2\backslash D$, called the artificial or source boundary, we easily see that any coefficient function $c:\Gamma\to \mathbb{C}^2$ would generate a smooth solution of \eqref{navier} in $D$ by the continuous superposition
\begin{align}\label{SL}
x\mapsto \int_{\Gamma}\Phi_{\varrho\omega^2}(x,y)c(y)\,\mathrm{d}s(y)=:\left(\Phi_{\varrho\omega^2}\ast_{|_{\Gamma}} c\right)(x)\ ,
\end{align}
which is the starting point for the MFS. We now refine
\begin{align}\label{trial}
\mathcal{H}(\omega):=\left\{u=\Phi_{\omega^2}\ast_{|_{\Gamma}} c_{u},\ v=\Phi_{\rho\omega^2}\ast_{|_{\Gamma}} c_{v}: (c_{u},c_{v})\in L^2(\Gamma)\times L^2(\Gamma)\right\}
\end{align}
from \eqref{ansatz} and prove in the following that the resulting approximation space is still sufficiently dense to finally recover exact eigenfunctions via boundary control. For this, we need some prerequisites concerning the construction of fundamental solutions for higher order PDEs. It also shows that the MFS ansatz in \eqref{SL} for approximating $u$ and $v$ separately is equivalent to recovering their difference $u-v$ via its synthesized fundamental solution. 
\begin{lem}
	If $\Phi_{\rho\omega^2}$ and $\Phi_{\omega^2}$ are fundamental solutions for the free space Navier equation \eqref{navier} with densities $\rho$ and $1$, respectively, then the function $\Phi_{\rho\omega^2,\omega^2}:=(\Phi_{\rho\omega^2}-\Phi_{\omega^2})/((1-\rho)\omega^2)$ is a fundamental solution for the fourth order operator $(\Delta^*+\omega^2)(\Delta^*+\rho\omega^2)$, where the operator product should be understood as composition.
	\begin{proof}
		Let $\phi \in C_c^{\infty}(D,\mathbb{C}^2)$ be an arbitrary bump function and set $\Phi_{\rho\omega^2,\omega^2}^y(x):=\Phi_{\rho\omega^2,\omega^2}(y-x)$ for some fixed $y\in\mathbb{R}^2$ (and likewise for its generating fundamental solutions). Then we can easily check
		\begin{align*}
		&\int_D\Phi_{\rho\omega^2,\omega^2}^y (\Delta^*+\omega^2)(\Delta^*+\rho\omega^2)\cdot\phi\,\mathrm{d}x \\
		=&\int_D \frac{\Phi_{\rho\omega^2}^y-\Phi_{\omega^2}^y}{(1-\rho)\omega^2} (\Delta^*+\omega^2)(\Delta^*+\rho\omega^2)\cdot\phi\,\mathrm{d}x \\
		=&\ \frac{1}{(1-\rho)\omega^2}\int_D\Phi_{\rho\omega^2}^y (\Delta^*+\rho\omega^2)(\Delta^*+\omega^2)\cdot\phi\,\mathrm{d}x \\
		&\ -\frac{1}{(1-\rho)\omega^2}\int_D\Phi_{\omega^2}^y (\Delta^*+\omega^2)(\Delta^*+\rho\omega^2)\cdot\phi\,\mathrm{d}x\\
		=&\ \frac{\Delta^*\phi(y)+\omega^2\phi(y)}{(1-\rho)\omega^2}-\frac{\Delta^*\phi(y)+\rho\omega^2\phi(y)}{(1-\rho)\omega^2} \\
		=&\ \phi(y)\ ,
		\end{align*}
		which proves the lemma.
	\end{proof}
\end{lem}
We are now ready to prove that for any eigenfunction pair $(u,v)$ of \eqref{ITP} with eigenfrequency $\omega$ 
lying in the first complex quadrant we can find approximations in $\mathcal{H}(\omega)$ with arbitrarily small boundary misfits. In particular, all the recovery criteria from Theorem \ref{theorem1} can be satisfied for corresponding ITEs. 

\begin{thm}\label{density}
	Let $w\in H^2(D,\mathbb{C}^2)$ be any distributional solution to the fourth order equation $(\Delta^*+\omega^2)(\Delta^*+\rho\omega^2)w=0$ for some frequency $\omega$ such that $0\leq \arg(\omega)<\pi/4$ (holds even more generally for $\mathrm{Im}\;\omega\geq 0$).
	Then there exists a sequence of elements $(u_m,v_m)_{m\in\mathbb{N}}\subset \mathcal{H}(\omega)$ from \eqref{trial} such that $(u_m-v_m)=:w_m\to w$ in $H^2(D,\mathbb{C}^2)$. If $\partial D$ is of class $C^{1,1}$, then in particular $\|w_m-w\|_{\sobolev^{\frac{3}{2}}(\partial D,\mathbb{C}^2)}\to 0$ and $\|\sigma(w_m-w)\nu\|_{\sobolev^{\frac{1}{2}}(\partial D,\mathbb{C}^2)}\to 0$.
	\begin{proof}
		Assume $\widetilde{w}\in \widetilde{H}^{-2}(D,\mathbb{C}^2)$, where the latter denotes the negative Sobolev space with compact support in $D$ identifying the dual space of $H^2(D,\mathbb{C}^2)$, is chosen such that 
		\begin{align}\label{duality}
		\int_D \widetilde{w}\cdot(\Phi_{\rho\omega^2,\omega^2}\ast_{|_{\Gamma}} c_1+\Phi_{\rho\omega^2}\ast_{|_{\Gamma}} c_2)\, \mathrm{d}x=0
		\end{align}
		for all $(c_1,c_2)\in L^2(\Gamma,\mathbb{C}^2)\times L^2(\Gamma,\mathbb{C}^2)$. The integral expression is hereby overloaded in notation with the corresponding duality pairing and by definition of $\Phi_{\rho\omega^2,\omega^2}$ from the previous lemma we see that the kernel of $\widetilde{w}$ contains functions of the form $(u-v)$ with $(u,v)\in\mathcal{H}(\omega)$. Therefore, if we can show that \eqref{duality} implies
		\begin{align}\label{check}
		\int_D\widetilde{w}\cdot w^*\mathrm{d}x =0
		\end{align}
		for every distributional solution $w^*\in H^2(D,\mathbb{C}^2)$ of $(\Delta^*+\omega^2)(\Delta^*+\rho\omega^2)w^*=0$, the Hahn-Banach theorem would yield the desired density claim since no other extension is possible. 
		
		For this, we define the auxiliary functions $w:=\Phi_{\rho\omega^2,\omega^2}\ast_{|_D}\widetilde{w}\in L^2(D,\mathbb{C}^2)\cap C^{\infty}(\mathbb{R}^2\backslash D,\mathbb{C}^2)$ and $v:=\Phi_{\rho\omega^2}\ast_{|_D}\widetilde{w}\in L^2(D,\mathbb{C}^2)\cap C^{\infty}(\mathbb{R}^2\backslash D,\mathbb{C}^2)$, where the gain of regularity in $D$ results from the fact that the convolution with $\Phi_{\rho\omega^2,\omega^2}$ or $\Phi_{\rho\omega^2}$ are pseudo-differential operators of order $-2$ due to their logarithmic singularity type. Rewriting \eqref{duality}, we obtain
		\begin{align*}
		\int_{\Gamma}c_1(y)\cdot w(y)+c_2(y)\cdot v(y)\,\mathrm{d}s(y)=0
		\end{align*}
		for all $(c_1,c_2)\in L^2(\Gamma,\mathbb{C}^2)\times L^2(\Gamma,\mathbb{C}^2)$ and setting one of the coefficient functions to zero, respectively, we may conclude that $w_{|\Gamma}=v_{|\Gamma}=0$. Using the pointwise estimate  $|\sqrt{r}(\partial_rv_p-\mathrm{i}(\rho\omega^2)v_p)(x)|\leq \sqrt{r}\|\partial_r(\Phi_{\rho\omega^2})_p-\mathrm{i}\rho\omega^2(\Phi_{\rho\omega^2})_p\|_{H^2(x-D,\mathbb{C}^{2\times2})} \|\widetilde{w}\|_{H^{-2}(D,\mathbb{C}^2)}$, and similarly for $v_s$, while employing standard differentiation properties and decay estimates for the resulting Hankel functions expansion within the first norm on the right, we deduce that Kupradze's radiation conditions are completely inherited by $v$ as convolution of the correspondingly radiating fundamental solution and some compactly supported distribution. By uniqueness of the exterior Navier problem for $\mathrm{Im}\;\omega\geq 0$ and $(\Delta^*+\rho\omega^2)v=\widetilde{w}$ in the sense of distributions, we may then conclude that $v=0$ outside of $\Gamma$. Due to analyticity, $v$ even needs to vanish completely in $\mathbb{R}^2\backslash \overline{D}$ because the right hand side of the Navier equation is zero here by assumption. Similarly, we want to prove in the following that $w\in H^2_0(D)$ 
		for justifying its role as a valid test function later:
		
		Using the definition of $\Phi_{\rho\omega^2,\omega^2}$, direct calculations show the distributional relations $(\Delta^*+\omega^2)w=v=(\Delta^*+\omega^2)(\Phi_{\omega^2}\ast_{|_D}v)$ as well as $(\Delta^*+\rho\omega^2)w=\Phi_{\omega^2}\ast_{|_D}\widetilde{w}=(\Delta^*+\rho\omega^2)(\Phi_{\omega^2}\ast_{|_D}v)$. This combines to $0=((\Delta^*+\rho\omega^2)-(\Delta^*+\omega^2))(w-\Phi_{\omega^2}\ast_{|_D}v)/(\rho-1)=w-\Phi_{\omega^2}\ast_{|_D}v$ and implies the additional representation $w=\Phi_{\omega^2}\ast_{|_D}v$. The same uniqueness and analyticity reasoning as for $v$ above but with frequency $\omega^2$ now yields $w=0$ in $\mathbb{R}^2\backslash \overline{D}$ again and by a bootstrap argument due to $v\in L^2(D,\mathbb{C}^2)$ we then conclude $w\in H^2_{loc}(\mathbb{R}^2,\mathbb{C}^2)$ or equivalently $w\in H^2_0(D)$. 
		Therefore we can find a sequence of bump functions $\{\phi_k\}_{k\in\mathbb{N}}$ 
		such that $(\Delta^*+\rho\omega^2)(\Delta^*+\omega^2)(w-\phi_k)\to 0$ in $H^{-2}(D,\mathbb{C}^2)$. Taking then any distributional solution $w^*\in H^2(D,\mathbb{C}^2)$ of $(\Delta^*+\omega^2)(\Delta^*+\rho\omega^2)w^*=0$ as mentioned in \eqref{check}, we may finally compute
		\begin{align*}
		\int_D \widetilde{w}\cdot w^*\,\mathrm{d}x &=\int_D \left((\Delta^*+\rho\omega^2)(\Delta^*+\omega^2)w\right)\cdot w^*\,\mathrm{d}x \\
		&=\lim_{k\to \infty} \int_D (\Delta^*+\rho\omega^2)(\Delta^*+\omega^2)\phi_k\cdot w^*\,\mathrm{d}x \\
		&=0\ ,
		\end{align*}
		where we especially incorporated in the first step that $(\Delta^*+\rho\omega^2)(\Delta^*+\omega^2)w=\widetilde{w}$ in the sense of distributions according to our original definition of $w$. Since $w^*$ was an arbitrary homogeneous solution, the desired density result for the interior domain is thereby proven. An application of standard trace theorems eventually takes the approximation result over to the boundary of $D$ in corresponding norms.
	\end{proof}
\end{thm}
In the next section we will focus on the numerical implications of these abstract findings and especially present its associated recipe of how to approximate ITEs in practice.

\subsection{The numerical setting}

The MFS aims to discretize \eqref{SL} in terms of a Riemann sum with unknown integration weights determined in some optimization process which we derive in the following. For arbitrary frequencies $\omega$ and some approximation order $m\in \mathbb{N}$ we focus on trial functions defined on $\overline{D}$ of the form
\begin{align*}
u_m(x)=\sum_{j=1}^{m}\Phi_{\omega^2}(x,y_j)\,(\widetilde{c}_u)_{1\leq i\leq 2,j}\ , \qquad 
v_m(x)=\sum_{j=1}^{m}\Phi_{\rho\omega^2}(x,y_j)\,(\widetilde{c}_v)_{1\leq i\leq 2,j}\ ,
\end{align*}
where $\{y_1,\dots,y_m\}\subset \Gamma$ 
are preselected source points and $\widetilde{c}_u,\widetilde{c}_v\in\mathbb{C}^{2\times m}$ are up to now unspecified block coefficient vectors.
Note $(u_m,v_m)$ are actually pseudo-elements of $\mathcal{H}(\omega)$ by using delta-distribution-type coefficient functions $c_u$ and $c_v$ along $\Gamma$, respectively. However, it can be readily seen that they are dense in $\mathcal{H}(\omega)$ for $m\to \infty$ in any interior Sobolev norm, cf. \cite{kleefeld2018method} where this was shown for the acoustic framework. In order to characterize the optimal pair $(\widetilde{c}_u,\widetilde{c}_v)\in\mathbb{C}^{2\times m}\times \mathbb{C}^{2\times m}$, 
we will now develop an MFS-based scheme that is supposed to meet the criteria listed in Theorem \ref{theorem1} from a numerical perspective. It uses the collocation concept as a discretized representative for the Sobolev-norm-conditioned quantities involved.

\begin{figure}
	\centering
	\includegraphics[scale=0.48]{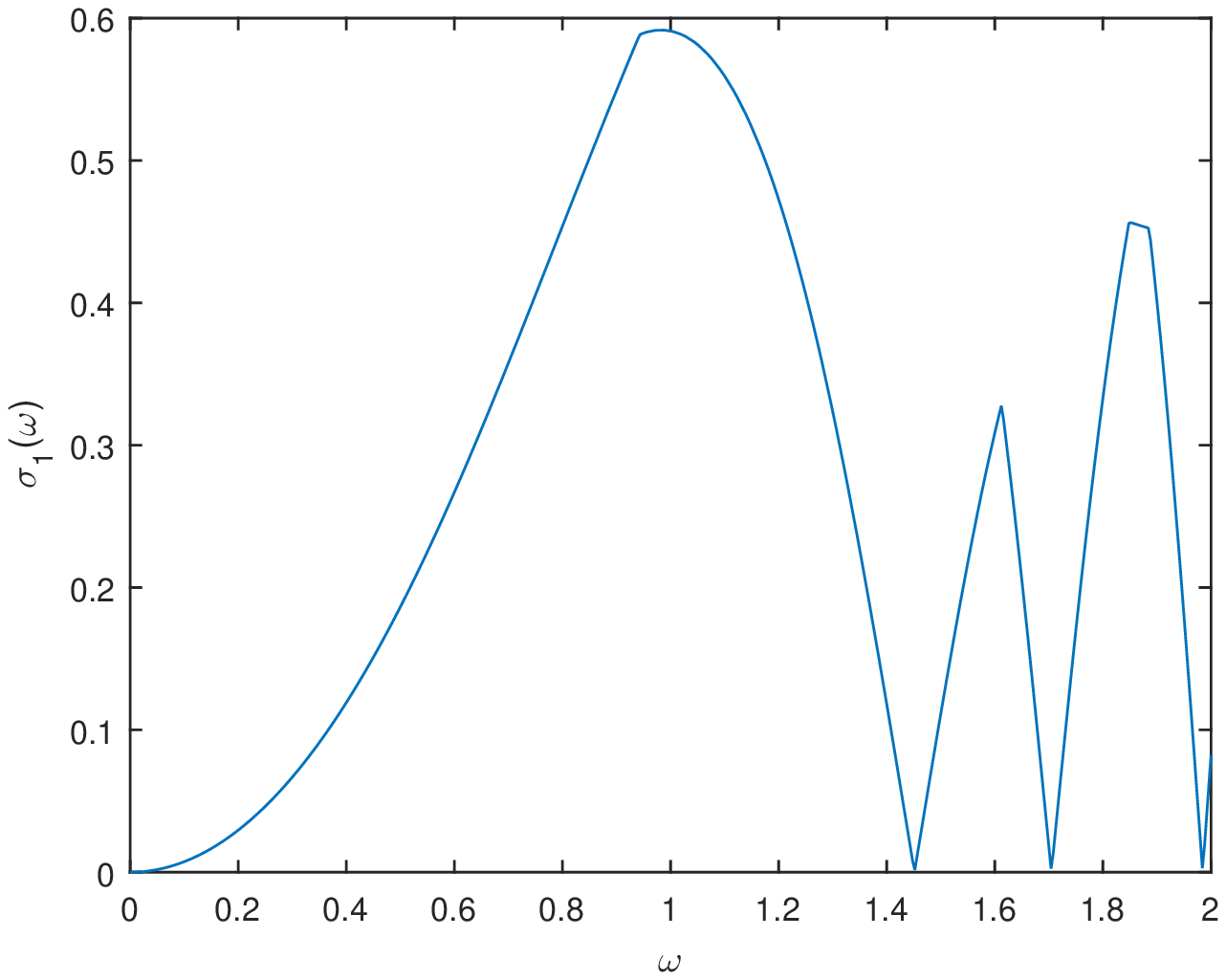}
	\includegraphics[scale=0.48]{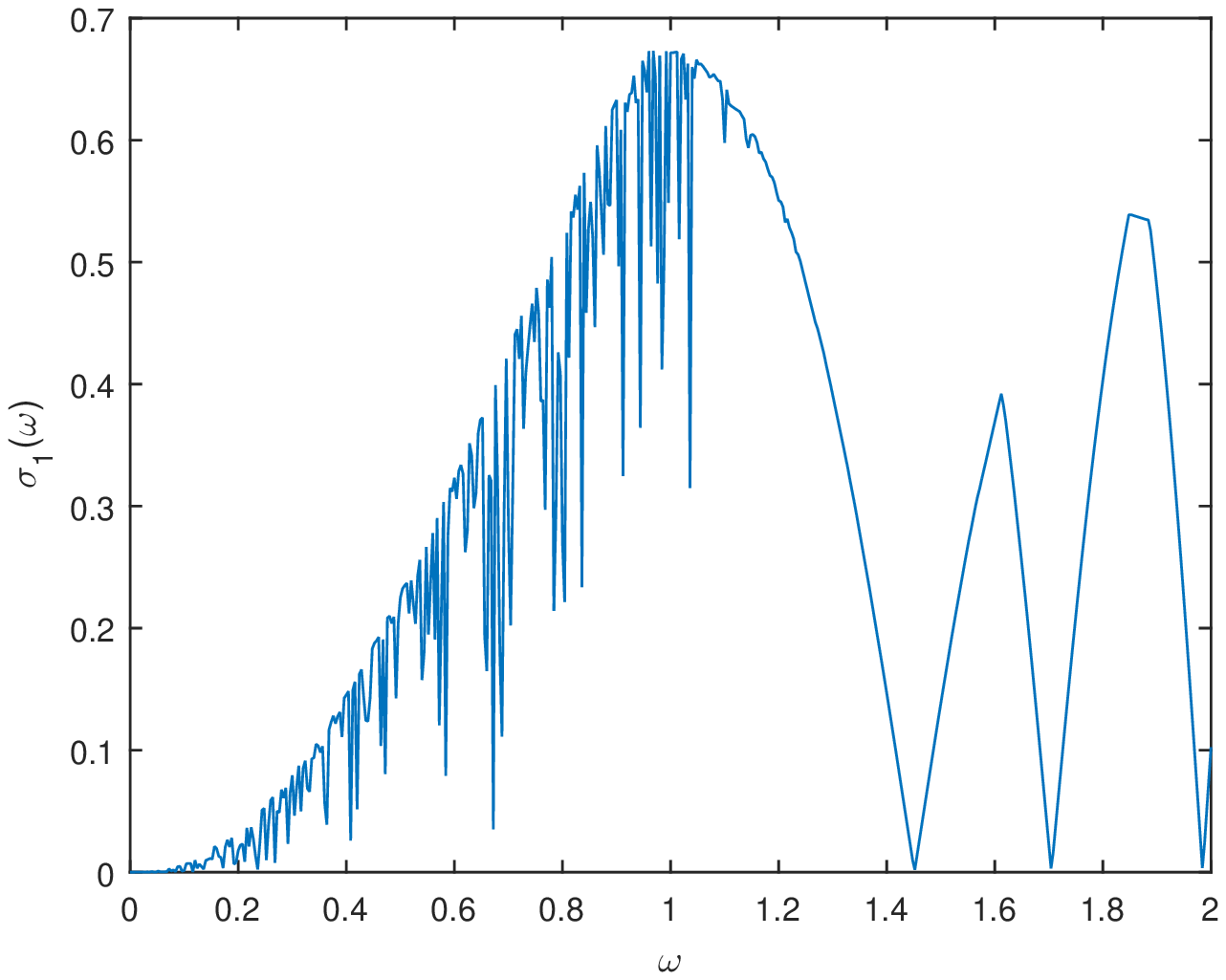}
	\caption{Exemplary plots of \eqref{svd} along the real axis for the disc of radius 0.5 as scatterer generated for $m=50$ (left) and $m=80$ (right). Both cases confirm 3 approximate ITEs $\omega_{50}$ within the interval $[0,2]$. However, for critically large $m$ the graph gets polluted by ill-conditioning artifacts.}\label{figure1}
\end{figure}

First, we treat the vanishing boundary misfit condition $(iii)$. Therefore we pick $m$ collocation points $\{x_1,\dots,x_m\}\subset \partial D$ and consider for fixed $\omega$ the constrained minimization of 
\begin{align*}
(\widetilde{c}_u,\widetilde{c}_v)\mapsto \sum_{i=1}^{m} |u_m(x_i)-v_m(x_i)|^2\ .
\end{align*}
This optimization must be subject to criterion $(ii)$ in a coherently discretized way for which we will take quadrature-like sample points $\{\widehat{x}_1,\dots,\widehat{x}_{m_I}\}\subset D$ with $m_I$ being large but fixed and demand
\begin{align*}
(\widetilde{c}_u,\widetilde{c}_v)\mapsto \sum_{i=1}^{m_I} |u_m(\widehat{x}_i)|^2+|v_m(\widehat{x}_i)|^2\approx 1\ .
\end{align*}
That is, we scale all eigenfunction candidate pairs to some almost constant but non-vanishing discretized interior norm to guarantee $(u_m,v_m)\ne 0$ for all $\omega$. Both sums above can be restated in compact matrix form, for which we define 
\begin{align}\label{blockext}
M(\omega):=
\begin{pmatrix}
\widetilde{\Phi}_{\omega^2} & \widetilde{\Phi}_{\rho\omega^2} \\
\sigma(\widetilde{\Phi}_{\omega^2})\nu  & \sigma(\widetilde{\Phi}_{\rho\omega^2})\nu  \\
\widehat{\Phi}_{\omega^2} & 0 \\
0 & \widehat{\Phi}_{\rho\omega^2}
\end{pmatrix}\in \mathbb{C}^{(4m+4m_I)\times 4m} \ ,
\end{align}
where the upper dense rows correspond to the boundary control matrices
\begin{align*} 
\left(\widetilde{\Phi}_{\varrho\omega^2}\right)_{2i-1,2j-1}&=\big(\Phi_{\varrho\omega^2}(x_i,y_j)\big)_{1,1}\ ,& \quad
\left(\widetilde{\Phi}_{\varrho\omega^2}\right)_{2i-1,2j}&=\big(\Phi_{\varrho\omega^2}(x_i,y_j)\big)_{1,2}\ , \\
\left(\widetilde{\Phi}_{\varrho\omega^2}\right)_{2i,2j-1}&=\big(\Phi_{\varrho\omega^2}(x_i,y_j)\big)_{2,1}\ ,& \quad
\left(\widetilde{\Phi}_{\varrho\omega^2}\right)_{2i,2j}&=\big(\Phi_{\varrho\omega^2}(x_i,y_j)\big)_{2,2}
\end{align*}
with $1 \leq i, j \leq m$ 
and co-normal derivatives affecting only the first arguments on the right hand side of involved definitions. A similar pattern applies for $\sigma(\widetilde{\Phi}_{\varrho\omega^2})\nu$. The diagonal lower part of \eqref{blockext} 
then similarly embodies the interior samples
\begin{align*}
\left(\widehat{\Phi}_{\varrho\omega^2}\right)_{2i-1,2j-1}&=\big(\Phi_{\varrho\omega^2}(\widehat{x}_i,y_j)\big)_{1,1}\ ,& \quad
\left(\widehat{\Phi}_{\varrho\omega^2}\right)_{2i-1,2j}&=\big(\Phi_{\varrho\omega^2}(\widehat{x}_i,y_j)\big)_{1,2}\ , \\
\left(\widehat{\Phi}_{\varrho\omega^2}\right)_{2i,2j-1}&=\big(\Phi_{\varrho\omega^2}(\widehat{x}_i,y_j)\big)_{2,1}\ ,& \quad
\left(\widehat{\Phi}_{\varrho\omega^2}\right)_{2i,2j}&=\big(\Phi_{\varrho\omega^2}(\widehat{x}_i,y_j)\big)_{2,2}\ ,
\end{align*}
where now $1\leq i \leq m_I$ and again $1 \leq j \leq m$. 
The benefit of this matrix reformulation is that we can now perform a $QR$-factorization 
\begin{align*}
M(\omega)=Q_M(\omega)R_M(\omega)=
\begin{pmatrix}
Q(\omega) \\ Q_I(\omega)
\end{pmatrix}
R_M(\omega)\ 
\end{align*}
with $Q\in \mathbb{C}^{4m\times 4m}$ and $Q_I\in \mathbb{C}^{4m_I\times 4m}$, which conveniently reflects the above minimal-boundary-to-maximal-interior coupling through the unitary property of $Q_M$. This observation goes back to Betcke and Trefethen, see \cite{betcke2005reviving}, who analyzed  Dirichlet eigenvalues with a slightly different ansatz but still via boundary control. Our common ingredient left is then to find local minima of
\begin{align}\label{svd}
\omega \mapsto \min_{r\in \mathbb{C}^{4m}, |r|=1}|Q(\omega)r|=:\sigma_1(\omega)\ .
\end{align}
Those solutions that are (almost) roots for the smallest singular value $\sigma_1(\omega)$, see Figure \ref{figure1} (left), will be denoted by $\omega_m$ to trace back to the underlying trial space dimension of $(u_m,v_m)$. Our described solution algorithm will be called modified MFS and its output, $\omega_m$, approximate ITEs. In the spirit of $(i)$ from Theorem \ref{theorem1} we hope them to converge for $m\to \infty$ also in practice and Lemma \ref{lemma} would then provide a measure for their convergence speed in terms of the continuous analogon of $\sigma_1(\omega_m)$. The next section will demonstrate this for some exemplary scatterers. 


\section{Numerical results}

We use the modified MFS from the previous section to compute some ITEs for a disc-, ellipse-, kite- and square-shaped scatterer $D$ whose (collocation) boundaries are given by
\begin{align*}
\partial D_d:=&\begin{pmatrix} 0.5 \cos(t) \\ 0.5 \sin(t) \end{pmatrix}\ , \quad t\in [0,2\pi)\ ,\\
\partial D_e:=&\begin{pmatrix} 0.5\cos(t) \\ \sin(t) \end{pmatrix}\ , \quad t\in [0,2\pi)\ , \\
\partial D_k:=&\begin{pmatrix} 0.75\cos(t)+0.3 \cos(2t) \\ \sin(t) \end{pmatrix}\ , \quad t\in [0,2\pi)\ , \\
\partial D_s:=&\ \partial\left([-0.5,0.5]\times[-0.5,0.5]\right)\ ,
\end{align*}
respectively. Note that $D_s$ is actually not covered by our ITE approximation analysis due to its corners, but still a very typical scatterer feasible for our algorithm since collocation is invisible with respect to the boundary regularity except for singular points that need to be excluded. Throughout we fix the constitutive parameters
\begin{align*}
\mu=\frac{1}{16}\ ,\quad \lambda=\frac{1}{4}\ ,\quad \rho=4
\end{align*} 
which were also used in the context of \cite{ji2018computation} and therefore serve as independent reference values for our exemplary MFS-based findings listed in Figure \ref{table:1}. Those approximate the first four real-valued ITEs from each scatterer and were individually obtained via the source point boundaries 
\begin{align*}
\Gamma_d:=&\begin{pmatrix} \cos(t) \\ \sin(t) \end{pmatrix}=2\cdot\partial D_d\ , \quad t\in [0,2\pi)\ , \\
\Gamma_e:=&\begin{pmatrix} 0.95\cos(t) \\ 1.9\sin(t) \end{pmatrix}=1.9\cdot\partial D_e\ , \quad t\in [0,2\pi)\ , \\
\Gamma_k:=&\begin{pmatrix} \left(1.2\cos(t)+0.48 \cos(2t)\right) \\ 1.6\sin(t) \end{pmatrix}=1.6\cdot\partial D_k\ , \quad t\in [0,2\pi)\ , \\
\Gamma_s:=&\partial\left([-0.65,0.65]\times[-0.65,0.65]\right)=1.3\cdot\partial D_s\ .
\end{align*}

\begin{figure}
	\centering
	\scalebox{0.8}{ 
		\begin{tabular}{ | l | l | l | l | l |}
			\hline
			\textbf{Shape} &  \textbf{ITE 1} & \textbf{ITE 2}& \textbf{ITE 3}& \textbf{ITE 4}\\ \hline
			disc					&
			1.451304027606383		&
			1.704638247023373		&
			1.984530256321993		&
			2.269112085458542		\\ \hline
			ellipse					&
			1.296728136516			&
			1.302785814026			&
			1.540896035208			&
			1.565151107263			\\ \hline
			kite					&
			0.947					&
			1.047					&
			1.111					&
			1.235					\\ \hline
			square					&
			1.3938					&
			1.6182					&
			1.8020					&
			1.9362					\\ \hline
		\end{tabular}
	}
	\caption{Approximations of the first four real-valued ITEs (counting without multiplicity) for some exemplary scatterers obtained via the modified MFS with material parameter $\mu=1/16,\ \lambda=1/4,\ \rho=4$: One clearly sees that the more advanced the shape of the boundary becomes, the less ITE digits can be effectively recovered.}
	\label{table:1}
\end{figure} 

The extracted scaling factors $\{2, 1.9, 1.6, 1.3\}$ were preselected to follow the conclusions from \cite{kleefeld2018method} and approach unity the more scattering shapes seem to deviate from the disc that is considered as the most promising state. Having thus set all the necessary computational contours for the modified MFS, both the $m$-dependent collocation points $\{x_1,\dots,x_m\}$ as well as the auxiliary sources $\{y_1,\dots,y_{2m}\}$ were distributed equiangular on their corresponding boundaries, e.g. equidistant with respect to $t$ if our representation above allows. However, it is well known that the individual optimization of $\Gamma$ itself and its source point distribution do have a noticeable impact on the MFS approximation quality as was shown in \cite{alves2009choice}, but resulting in a more advanced non-linear problem for each scatterer in total. In contrast, concerning the interior points, $m_I=10$ of them were fixed randomly in a centered disc with radius $0.5$ as their overall location and number turns out not to affect the ITE output significantly.

\begin{figure}
	\centering
	\includegraphics[scale=0.8]{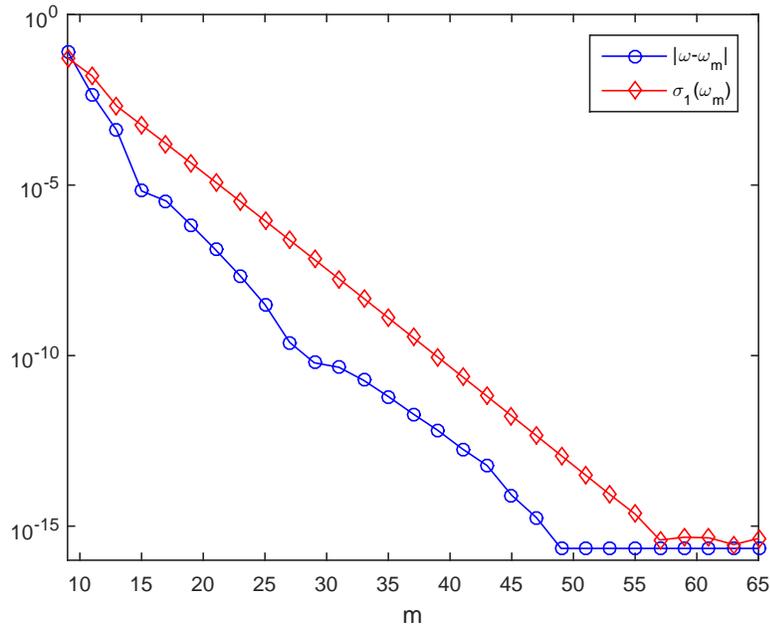}
	\caption{Exponential decay of the modified MFS with respect to the absolute ITE deviation from 1.451304027606383 being the smallest real-valued one of the disc with radius 0.5 as scatterer. A similar behavior for the boundary data misfit of the corresponding approximate eigenfunctions in terms of the smallest singular value is shown. Their linear correlation in the plot confirms the continuous analogon from Lemma \ref{lemma}.} \label{figure2} 
\end{figure}

With these input arrangements and the focus on real-valued ITEs first, the modified MFS can be successfully exploited by incrementing $m$ within $40\leq m\leq 80$, where the lower bound just provided an averaged value for when to expect good results. However, exceeding the given $D$-specific threshold for our setup, the graph of \eqref{svd} always began to suffer drastically from the intrinsic ill-conditioning effects of the discretization matrix $M(\omega)$ in \eqref{blockext} via impeding oscillations and finally lead to unreliable approximations apart, see Figure \ref{figure1} (right). Obtained via the regime in between, we believe our given cut-off results from Figure \ref{table:1} to be correct up to that point modulo round-off errors compared to the exact ITE mantissa since they correspond to the nonfluctuating digits within the modified MFS output $\omega_m$ when $m\nearrow 80$. Some exact reference values for $D_d$ corresponding to rotational-symmetric eigenfunctions can be obtained by computing roots of \eqref{det} which even confirms our approximation for its smallest ITE to be exact up to machine precision, see Figure \ref{figure2} for the convergence history. In general, the recoverable accuracy strongly correlates with the scattering shape and especially with the existence of corners. These technical observations mostly agree with the ones from the acoustic case analyzed in \cite{kleefeld2018method,kleefeld2018computing} and only differ in the retarded yet $D$-specific convergence regime in $m$ for starting the actual ITE approximations which was about $m\geq 20$ before. It manifests the fact that the eigenfunctions from elasticity are vector-valued and thus numerically slightly more expensive. 

\begin{figure}[!h]
	\centering
	\includegraphics[scale=0.47]{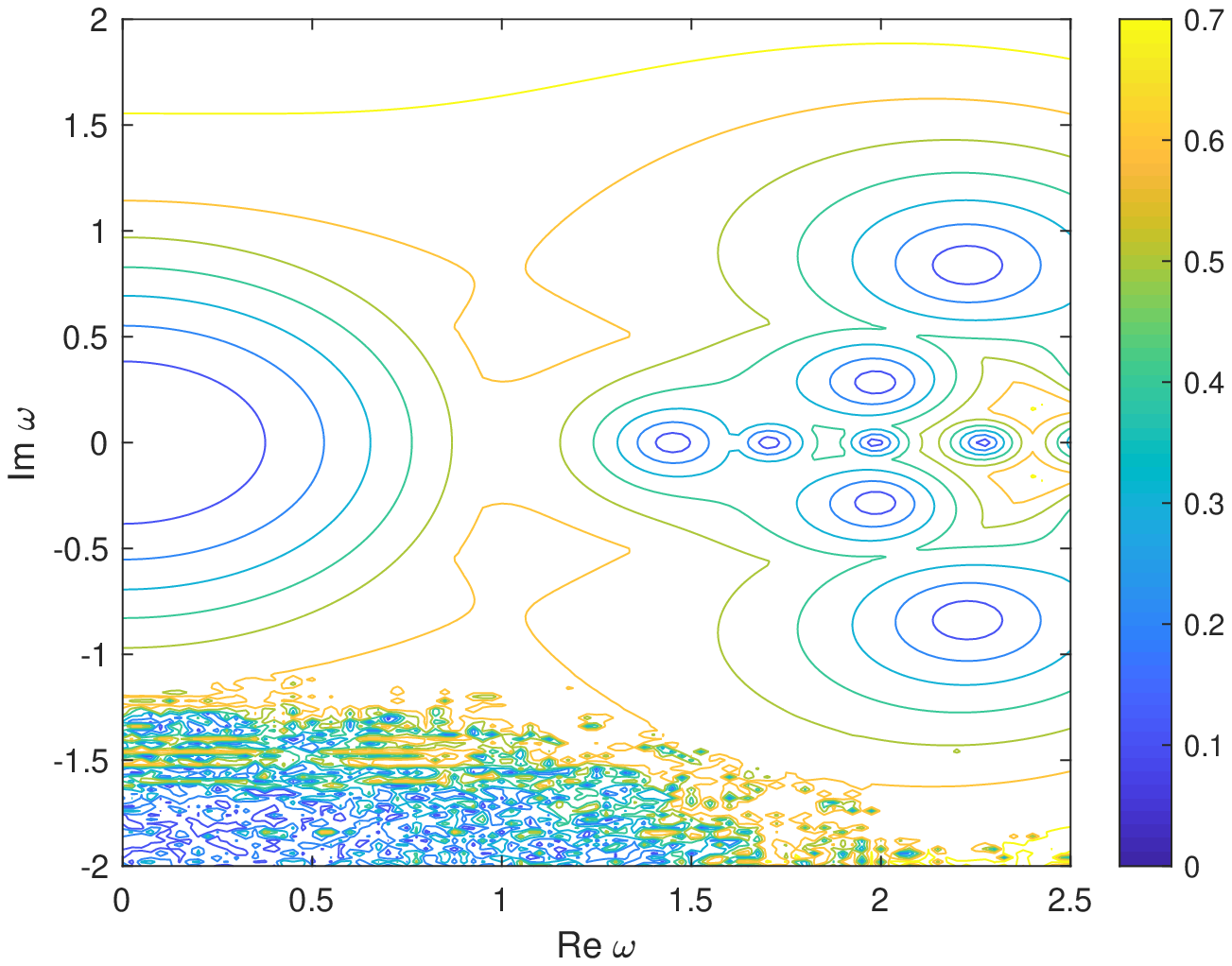}
	\includegraphics[scale=0.47]{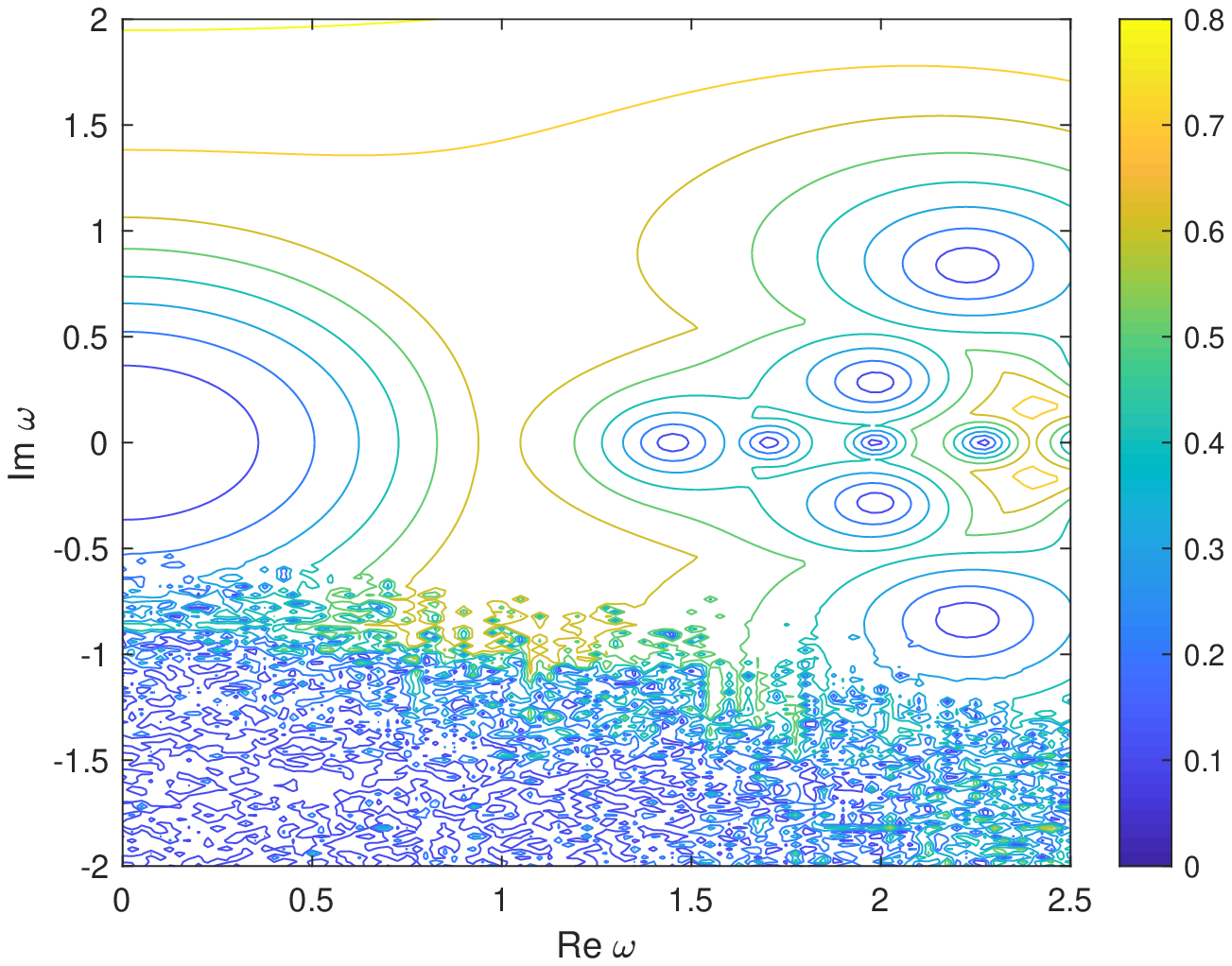}
	\includegraphics[scale=0.47]{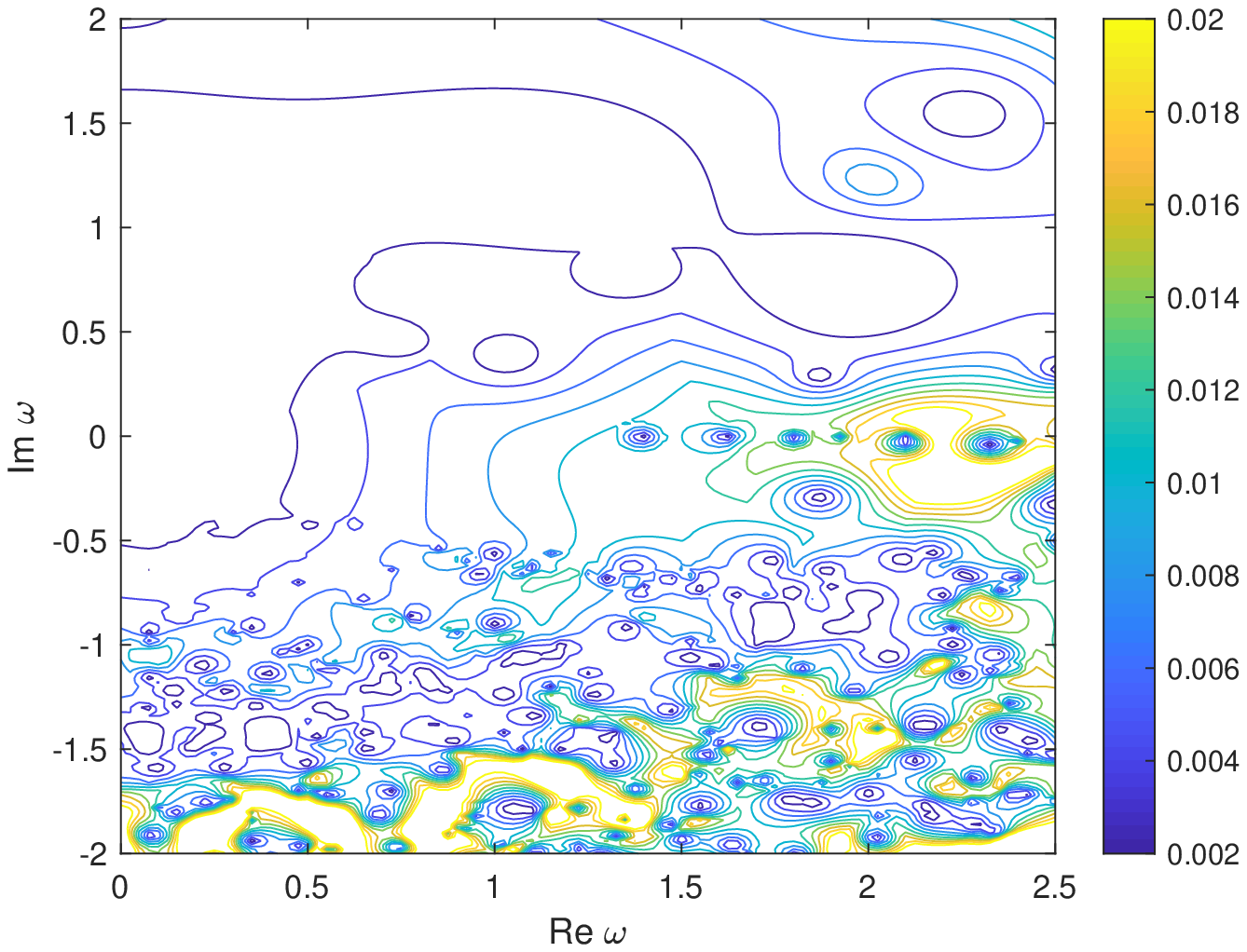}
	\includegraphics[scale=0.47]{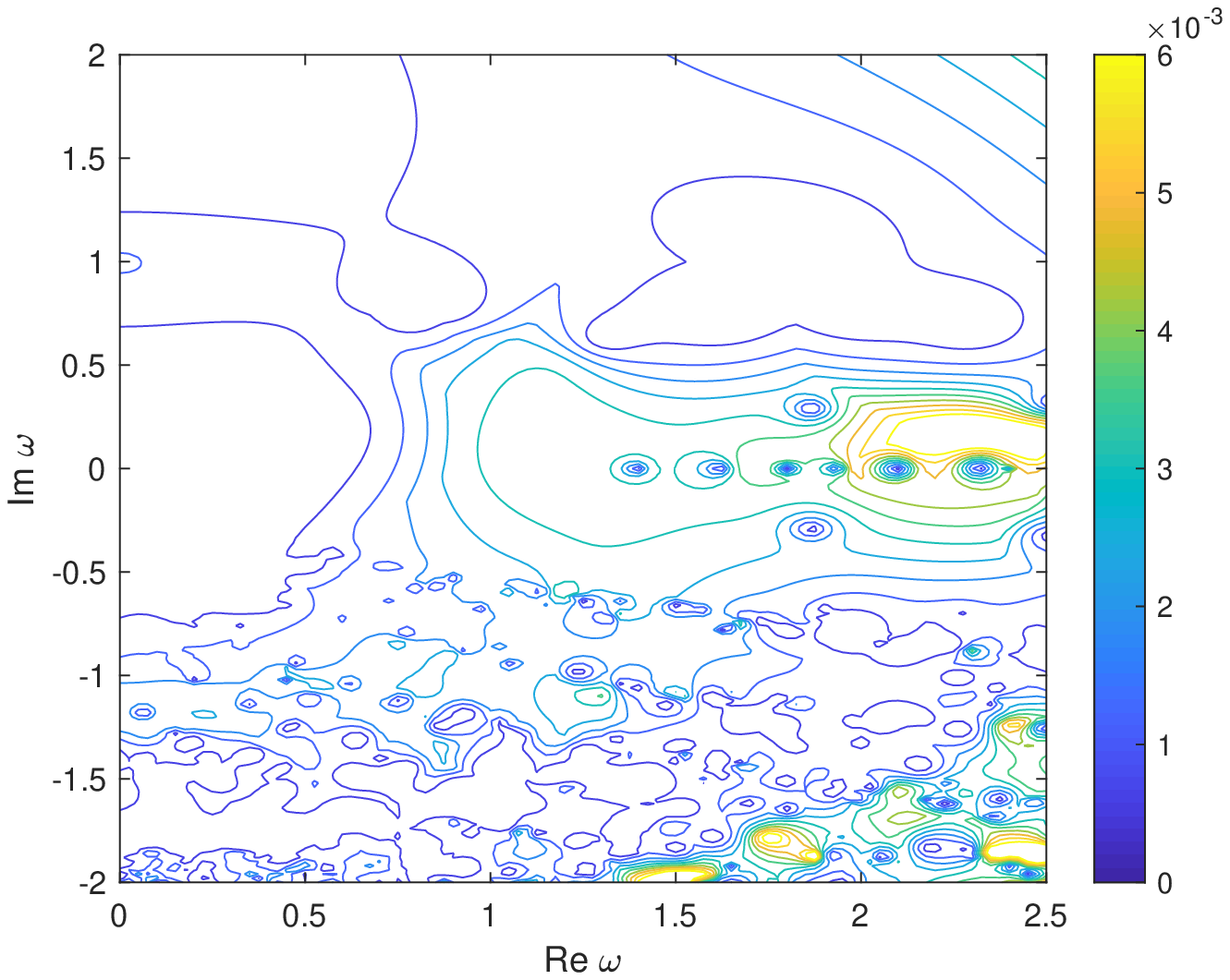}
	\caption{Contour extract from the complex plane evaluating \eqref{svd} for the unit square (below) and the disc of radius 0.5 (above) as scatterers with $m=45$ (left) and $m=55$ (right), respectively. Approximate ITEs are denoted by the centers of emerging concentric circles and spoilt primarily from the lower half space through the ill-conditioning artifacts already encountered in Figure \ref{figure1}.} \label{figure3} 
\end{figure}

As in the particular scope of this paper, we also want to present some approximations from the complex-valued eigenvalue spectrum. Therefore Figure \ref{figure3} plots \eqref{svd} in the complex range  $0\leq \mathrm{Re}\,\omega\leq 2.5,\ -2\leq\mathrm{Im}\,\omega\leq 2$ for the disc $D_d$ and the square $D_s$ since only their distribution of ITEs seems not that dense among our analyzed scatterers according to Figure \ref{table:1} and is thus better suited for broad contour visualizations. We can immediately see that ITEs are correctly displayed in conjugated pairs and as an example the closest of them with respect to the positive imaginary axis are computed for the two scatterers to be $1.8624+0.3104\mathrm{i}$ and $1.987178187576699+0.283125784408650\mathrm{i}$, respectively. The correctness of all given ITE digits for $D_d$ can even be confirmed again by comparing with the correspondingly coinciding root of \eqref{det}. Regarding accuracy, the modified MFS does not show any remarkable difference for its extension to the upper half space of the complex plane. However, although the plot indicates a certain axial symmetry with below, one clearly sees that inaccuracies due to large condition numbers of the underlying matrix \eqref{blockext} dominantly propagate from the lower half space upwards for increasing $m$, giving a first hint where the oscillations in Figure \ref{figure1} (right) originate from. Fortunately, ITEs always arise in conjugated pairs and because of Theorem \ref{density} we even know that restricting to complex numbers with non-negative imaginary part within our investigations, as already implemented in our initial ansatz \eqref{ansatz} and \eqref{trial}, is not limiting at all. 

To get a rough impression about how competitive our method is compared to others, we use the approximation results for real-valued $k^2$ from the paper \cite{ji2018computation} in which the authors also analyzed the scatterers $D_d$ and $D_s$ but via special finite element methods. With a mesh size of $h\leq 0.0125$, they obtained the approximations $k_d^2\approx 2.110723$ and $k_s^2\approx 1.9428781$ for the smallest ITE squared of the disc and the square, respectively. Allowing for possible round-off deviations in the last digit within our corresponding cut-off values, we obtain the tolerance ranges $k_d^2\in [2.106283380546506, 2.106283380546512]$ and $k_s^2\in [1.9423, 1.9429]$ by setting at most $m=80$ which may be converted via a boundary partitioning to a collocation point distance of the same order as the alternative mesh size. While the results regarding $D_s$ are of the same order for both solution approaches, the modified MFS clearly dominates in accuracy for $D_d$ and is therefore believed to do so for slightly perturbed scatterers $D$, too. 
\section{Conclusion}
We analyzed the method of fundamental solution in a stabilized version for the computation of complex-valued elastic transmission eigenvalues in two dimensions. Our theoretical studies show that the short algorithm shall produce spurious-free results in the limit whose approximation error per step we could quantify in terms of some discretized residual output value. Our numerical experiments confirm these expectations in practice within some scatterer-specific collocation point regime for which the method's feasibility is unaffected from ill-conditioning pollutions. This was tested for a collection of scatterers whose boundaries are easily-parametrizable. In accordance with the conclusions from previous works that analyzed our algorithm in the context of related transmission problems but restricted to real-valued eigenvalues so far, the best results, including also the complex spectrum from now, can still be obtained for the unit disc. The more the actual scattering shape then deviates from the disc, the less accurate approximations are finally achievable while even more collocation points are needed in total. However, depending on the choice of fundamental solution for generating the trial functions, different areas of the complex spectrum (in our case the lower half space) should be avoided as initial guess input for the algorithm due to ill-conditioning effects which are, however, not that restrictive because all eigenvalues arise in conjugated pairs.  
Generally, since the source boundary necessary for the MFS setup was individually preselected by intuition for simplicity, even more promising approximations can be achieved by its optimization which was, however, not covered in the scope of this paper. As a conclusion, the modified MFS is especially effective for quite regular domains and dominates here over many competitive methods in the general context of transmission eigenvalue problems.

\section*{References}	 
\bibliographystyle{abbrv}
\bibliography{Bibliography}

\end{document}